\documentclass[12pt]{article}
\usepackage{amssymb, amsmath}
\begin{document}
\voffset=0.0truein \hoffset=-0.5truein \setlength{\textwidth}{6.0in}
\setlength{\textheight}{8.8in} \setlength{\topmargin}{-0.2in}
\renewcommand{\theequation}{\arabic{section}.\arabic{equation}}
\newtheorem{thm}{Theorem}[section]
\newtheorem{definition}{Definition}[section]
\newtheorem{lemma}{Lemma}[section]
\newtheorem{pro}{Proposition}[section]
\newtheorem{cor}{Corollary}[section]
\newcommand{\n}{\nonumber}
\renewcommand{\o}{\omega}
\renewcommand{\O}{\Omega}
\newcommand{\e}{\varepsilon}
\newcommand{\vare}{\varepsilon}
\newcommand{\bb}{\begin{equation}}
\newcommand{\ee}{\end{equation}}
\newcommand{\bq}{\begin{eqnarray}}
\newcommand{\eq}{\end{eqnarray}}
\newcommand{\bqn}{\begin{eqnarray*}}
\newcommand{\eqn}{\end{eqnarray*}}
\title{On the apparition of singularities of vector fields
transported by volume preserving diffeomorphisms}
\author{Dongho Chae\thanks{This work was supported partially by  KRF Grant(MOEHRD, Basic
Research Promotion Fund) and the KOSEF Grant no.
R01-2005-000-10077-0.\newline
 Key Words: singularities,
 diffeomorphism, inviscid flows}\\
Department of Mathematics\\
              Sungkyunkwan University\\
               Suwon 440-746, Korea\\
              {\it e-mail : chae@skku.edu}}
 \date{}
\maketitle
\begin{abstract}
 We consider  possible generation of singularities of a vector field
 transported by diffeomorphisms with derivatives of uniformly bounded determinants.
 A particular case  of volume preserving diffeomrphism is the most
 important, since it has direct applications to the incompressible, inviscid
 hydrodynamics.
 We find
 relations between the directions of the vector field and the eigenvectors
 of the derivative of the back-to-label map near the singularity.
 We also find an invariant when we follow the motion of the integral curves of the vector
 field. For the 3D  Euler equations these results have
 immediate implications about the directions of the vortex
 stretching and the material stretching near the possible singularities.
 We also have similar applications to the other inviscid, incompressible fluid
 equations such as the 2D quasi-geostrophic equation and
 the 3D magnetohydrodynamics  equations.

\end{abstract}

\section{Vector fields transported by diffeomorphisms}
 \setcounter{equation}{0}

\subsection{Statement of the theorems}

 Let $D$
be  a domain in $\Bbb R^n$, and $T\in (0, \infty]$. Suppose that for
all $t\in [0, T)$ the mapping $a\to X(a,t)$ is a diffeomorphism on
$D$. We denote by $A(\cdot ,t)$ the inverse mapping of $X(\cdot
,t)$, satisfying
$$ A(X(a,t),t)=a, \quad X(A(x,t),t)=x \quad \forall a,x \in D,
\quad\forall t\in [0,T).
$$
In the applications to hydrodynamics in the next section the
 mapping $\{ X(\cdot,t)\}$ is defined by a smooth velocity field $v(x,t)$
 through the system of the ordinary
differential equations:
 \bb\label{1.1} \frac{\partial X(a,t)}{\partial t }
=v(X(a,t),t) \quad ;\quad X(a,0)=a\in D\subset \Bbb R^n.
 \ee
In such a case we say the `particle trajectory' map $X(\cdot, t)$
and its inverse, the `back-to-label' map $A(\cdot,t)$  are generated
by the fluid velocity field $v(x,t)$.

\begin{definition}
We call a parametrized vector field $W(\cdot, \cdot):D\times
[0,T)\to\Bbb R^n$ is transported by  a differentiable mapping $
X(\cdot ,t)$ from $D$ into itself for all $t\in [0,T)$ if
 \bb\label{1.2}
W(X(a,t),t)=\nabla_a X(a,t) W_0(a)
 \ee
 holds  for all $(a,t)\in
D\times [0, T)$, where we set $W_0 (x)=W(x,0)$.
\end{definition}

We note that (\ref{1.2}) corresponds to the  well-known  vorticity
transport formula for the incompressible Euler equations. Actually
it is well-known(see e.g. \cite{maj2}) that (\ref{1.2}) is
equivalent to saying that the vector field $W(x,t)$ satisfies the
system of differential equations:
 \bb\label{1.2a}
\left\{ \aligned &\frac{\partial W}{\partial t} +(v\cdot \nabla )W
=(W\cdot\nabla )v,\\
&W(x,0)=W_0 (x)
\endaligned \right.
\ee on $D\times [0, T)$, where $v(x,t)$ is defined from $X(\cdot
,t)$
 by (\ref{1.1}).
In this paper we are concerned on the study of the direction of
$W(x,t)$ and the directions of stretching/compressions induced by a
`generalized volume preserving' mapping $X(\cdot ,t)$  near possible
singularities  in (\ref{1.2})(or, equivalently in (\ref{1.2a})).
This will be done efficiently in terms of the derivative of its
inverse mapping $A(\cdot,t)=X^{-1}(\cdot,t)$. The main motivation of
the current study is to understand the dynamic relation between the
vortex stretching and the material stretching when we approach to
possible singularities in the 3D incompressible Euler equations and
other inviscid flows.

\begin{thm}
Let $W(x,t)$ be a vector field on $D\subset \Bbb R^n $ defined for
$t\in [0, T)$. We set $W_0(x)=W(x,0)$ with $\|W_0\|_{L^\infty (D)}
<\infty$. Suppose $W(x,t)$ is transported by a volume preserving
diffeomorphism $\{X(\cdot ,t)\}_{t\in [0, T)}$ on $D\subset \Bbb
R^n$, whose inverse is $A(\cdot ,t)$. We assume that
 \bb\label{hyp}
\sup_{(a,t)\in D\times [0,T)}|\mathrm{det}(\nabla_a
  X(a,t))|<\infty.
 \ee
Let us set by $\{ (\lambda_j ,e_j)\}_{j=1}^n$  the eigenvalue and
eigenvector pairs of the symmetric, positive definite matrix
$$M(x,t)= (\nabla A(x,t))^* \nabla A(x,t)
$$
with the order of magnitude
 \bb\label{1.13}
   \lambda_1 \geq \lambda_2  \geq
\cdots \geq \lambda_n >0.
 \ee
  Suppose there exists a sequence
$(x_k,t_k)$ and $(\bar{x} ,\bar{t} )$ in $\bar{D}\times [0,T]$ such
that $\lim_{k\to \infty} (x_k,t_k)=(\bar{x},\bar{t})$ and
$$
\lim_{k\to \infty}|W(x_k,t_k)|=\infty.
$$
Then, necessarily
  \bb\label{1.3}
 \lim_{k\to \infty} \lambda_1 (x_k,
t_k)=\infty, \quad \mbox{and} \quad \lim_{k\to \infty} \lambda_n
(x_k, t_k)=0.
 \ee
Let $m$ be the largest number in $\{ 1, \cdots, n\}$ such that
 \bb\label{1.3a}
 \lim_{k\to \infty}
\lambda_j (x_k,t_k)>0 \quad \forall j\in \{1, \cdots, m\}.
 \ee
We set by $\Xi (x,t)$ the direction field of $W(x,t)$ defined by
$$\Xi (x,t)=\frac{W(x,t)}{|W(x,t)|}\quad\mbox{whenever}\quad  |W(x,t)|\neq 0.
$$
 Then,
\bb\label{1.5}
 \lim_{k\to \infty } e_j (x_k,t_k)\cdot \Xi (x_k,t_k)=0 \quad
 \forall j\in \{1, \cdots, m\}.
\ee
 \end{thm}
 \noindent{ \textsc{Remark} 1.1} Let $v\in \Bbb R^n$. Then, the quantity $|\nabla_a X(a,t)
 v|/|v|$ has the meaning of the rate of stretching(compression) in the direction of $v$
  induced by the trajectory mapping
 $X(\cdot, t)$ if the quantity is bigger(less) than 1. Indeed, let
 $\{\gamma_0 (s)\}_{s\in (-\vare ,\vare )}$ be a curve in $\Bbb R^n$ such that
 $$\gamma_0 (0)=a, \quad \left. \frac{\partial \gamma_0 (s)}{\partial s}\right|_{s=0}=v .$$
 We set $X(\gamma_0 (s),t)=\gamma(s,t)$. Then,
 \bb
 \frac{\partial
 \gamma(s,t)}{\partial s}=\nabla_a X(\gamma_0 (s),t)\frac{\partial \gamma_0
 (s)}{\partial s},
 \ee
 and
 \bb
 \frac{|\nabla_a X(a,t) v|}{|v|}=\left. \frac{\left|\nabla_a X(\gamma_0 (s),t)\frac{\partial \gamma_0
 (s)}{\partial s}\right|}{\left|\frac{\partial \gamma_0
 (s)}{\partial s}\right|}\right|_{s=0}
=\left.\frac{\left|\frac{\partial \gamma(s,t)}{\partial
s}\right|}{\left|\frac{\partial \gamma_0
 (s)}{\partial s}\right|}\right|_{s=0}
 = \left|\frac{\partial \gamma (s,t)}{\partial \gamma_0
 (s)}\right|_{s=0},
 \ee
 which provides us with the desired interpretation.
 By the Rayleigh-Ritz theorem(\cite{hor}) and the fact that $X(\cdot, t)$ is a
 diffeomorphism we have
 \bqn
 \lambda_1(x,t)&=&\max\{ \lambda_1 (x,t), \cdots, \lambda_n
 (x,t)\}\n \\
 &=&\sup_{v\neq 0} \frac{v^* M(x,t) v}{|v|^2}
 =\sup_{v\neq 0}
 \frac{|\nabla A(x,t) v|^2}{|v|^2}\\
 &=&\sup_{w\neq
 0}\frac{|w|^2}{|\nabla_a X(a,t) w|^2}
 =\frac{1}{\inf_{w\neq 0}
 \frac{|\nabla_a X(a,t) w|^2}{|w|^2}}.
 \eqn
 Hence,
\bb \inf_{v\neq 0}
 \frac{|\nabla_a X(a,t) v|}{|v|}=\frac{1}{\sqrt{\lambda_1
 (x,t)}},
 \ee
where $x=X(a,t)$.
 Similarly, for $\lambda_n (x,t)=\min\{ \lambda_1 (x,t), \cdots,
 \lambda_n (x,t)\}$, we obtain
 \bb
\sup_{v\neq 0}
 \frac{|\nabla_a X(a,t) v|}{|v|}=\frac{1}{\sqrt{\lambda_n
 (x,t)}}
 \ee
 with $x=X(a,t)$. In particular, in the case of  det$(\nabla_a X(a,t))\equiv
 1$(incompressible flow), the quantity $1/\sqrt{\lambda_1 (x,t)} (<1)$  has the meaning of the
maximum
 compression rate,
 while $1/\sqrt{\lambda_n (x,t)} (>1)$ has the meaning of the maximum stretching rate at
 $(x,t)$, except the case $\nabla_a X(a,t)=I$, where $a=A(x,t)$.\\
\ \\
\noindent{ \textsc{Remark} 1.2} In particular (\ref{1.3}) implies
that the directions of  the infinite stretching rate and the zero
compression rate are mutually
orthogonal to each other.\\
\ \\
\noindent{ \textsc{Remark} 1.3} Since $\lim_{k\to
\infty}\lambda_j(x_k,t_k) =0$ for all $j=m+1,
 \cdots, n$ by the hypothesis of the above theorem, the conclusion
 (\ref{1.5}) implies that as $(x_k,t_k) \to (\bar{x}, \bar{t})$ the sequence of direction
 vectors
 $\{ \Xi(x_k,t_k)\}$ tends to be on the
 hypersurface spanned by the vectors with the directions of infinite
 stretching rates.\\
 \ \\
The first part of the following theorem could be regarded as a
generalization of the well-known Helmholtz vortex theorem for the
incompressible Euler equations(\cite{hlm}).
\begin{thm}
Suppose $W(x,t)$ is a vector field transported by a diffeomorphism
$\{ X(\cdot, t)\}$, $t\in [0, T)$. Let $\{ \gamma_0(s)\}_{s\in I}$
be an integral curve of $W(x,0)$, then $\gamma (s,t):=X(\gamma_0
(s),t)$ is also an integral curve of $W(\gamma(s,t),t)$. Moreover,
we have the following invariant:
 \bb\label{1.8}
 \frac{|W (\gamma (s,t),t)|}{ \left|\frac{\partial \gamma (s,t)}{\partial s}
 \right|}= \frac{|W_0 (\gamma_0(s))|}{
 \left|\frac{\partial \gamma_0(s)}{\partial s}
 \right|}.
\ee
\end{thm}
\noindent{ \textsc{Remark} 1.4} We will see in the proof of the
above theorem that the invariant (\ref{1.8}) is due to the  fact
that the integral curve of of a vector field has re-parametrization
symmetry.

\subsection{Proof of the main theorems}

{\bf Proof of Theorem 1.1} The vector field transport formula
 $$ W (X(a,t),t)=\nabla_a X(a,t)W_0 (a),
 $$
can be written as
 \bb\label{1}
 \nabla A(x,t)W(x,t)= W_0 (A(x, t))
 \ee
 in terms of $A(x,t)=X^{-1}(x,t)$. Hence,
\bq\label{1.12}
 |W_0 (A(x,
 t))|^2&=& W(x,t)^* (\nabla A(x,t))^* \nabla A(x,t)W(x,t)\n \\
 &=& |W(x,t)|^2\Xi(x,t)^* M(x,t)\Xi (x,t)\n \\
 &=& |W(x,t)|^2\left( \lambda_1 (x,t)\widetilde{\Xi}_1^2(x,t)+\cdots +\lambda_n (x,t)
 \widetilde{\Xi}_n^2(x,t)\right),\n \\
\eq
  where we set
  $$\widetilde{ \Xi }(x,t)=O(x,t)\Xi (x,t), \quad O^* M
  O=\mathrm{diag}(\lambda_1 ,\cdots, \lambda_n ).
  $$
 Namely, the $n\times n$ orthogonal matrix $O(x,t)$ diagonalizes the positive
 definite, symmetric matrix $M$.
 By the hypothesis (\ref{hyp}),
  \bq\label{1.14}
  \Lambda:&=&\inf_{(x,t)\in D\times [0, T)}\left[
  \lambda_1(x,t)\cdots  \lambda_n(x,t)\right]\n \\
  &=&\inf_{(x,t)\in D\times [0, T)}\mathrm{det}\,M
  =
  \inf_{(x,t)\in D\times [0, T)}\mathrm{det}(\nabla
 A(x,t))^2\n \\
 &=&\frac{1}{\sup_{(x,t)\in D\times [0, T)}\mathrm{det}(\nabla_a
 X(a,t))\big|^2_{a=A(x,t)}}>0.
 \eq
 By definition
  \bb\label{1.15}
\widetilde{\Xi}_1 ^2 (x,t)+\cdots +\widetilde{\Xi}_n ^2 (x,t)=|O
(x,t) \Xi (x,t)|^2= |\Xi (x,t)|^2=1.
 \ee
 Hence, from (\ref{1.12}) and the inequality $\frac{a_1+\cdots +a_n}{n}\geq
 (a_1\cdots a_n )^{\frac1n}$ for $a_1,\cdots, a_n \geq 0$, we obtain
 that
 \bq\label{1.16}
 \frac{\|W_0 \|_{L^\infty}^2}{|W(x,t)|^2}&\geq&
\lambda_1(x,t)\widetilde{\Xi}_1^2(x,t)+\cdots +\lambda_n(x,t)\widetilde{\Xi}_n^2(x,t)\n \\
 &\geq & n\left( \lambda_1
 (x,t)\widetilde{\Xi}_1^2(x,t)  \cdots \lambda_n
 (x,t)\widetilde{\Xi}_n^2(x,t)\right)^{\frac1n}\n \\
 &=&n\left( \lambda_1 (x,t)\cdots \lambda_n (x,t)\right)^{\frac1n}
 \left(\widetilde{\Xi}_1^2(x,t)  \cdots \widetilde{\Xi}_n^2(x,t)\right)^{\frac1n}\n \\
 &\geq&n \Lambda^{\frac1n}\left|\widetilde{\Xi}_1(x,t) \cdots \widetilde{\Xi}_n
 (x,t)\right|^{\frac2n}.
\eq
 Let $\{(x_k, t_k)\}$ be a sequence such that
  $(x_k ,t_k)\to (\bar{x} ,\bar{t})$ as $k \to \infty$, and
  $$ \lim_{k\to \infty}|W (x_k, t_k )|=\infty.$$
  Then, the first inequality of (\ref{1.16}) implies that
 \bb\label{1.17}
 \lim_{k\to \infty} \lambda_j (x_k,t_k )\widetilde{\Xi}^2_j
  (x_k,t _k )=0 \quad \forall j =1,\cdots , n.
  \ee
From (\ref{1.15}) we find that there exists $j_0 \in \{ 1, \cdots,
n\}$
 such that
 \bb\label{1.18}
 \lim\inf_{k\to \infty} |\widetilde{\Xi }_{j_0}| (x_k,t_k)\geq
\frac{1}{\sqrt{n}}.
 \ee
Hence,  from (\ref{1.13}) we have
 \bb
\lim_{k\to \infty}\lambda_{n} (x_k,t_k)\leq \lim_{k\to \infty}
\lambda_{j_0} (x_k,t_k)=0,
 \ee
which, in turn, implies by (\ref{1.14})  and (\ref{1.5}) that
 \bb
\lim_{k\to \infty}\lambda_{1} (x_k,t_k)=\infty.
 \ee
Thus we find there exists $m\in \{ 1, \cdots ,n-1\}$ satisfying
(\ref{1.3a}). Now  (\ref{1.17})  and (\ref{1.3a}) imply that
$$
\lim_{k\to \infty} \widetilde{\Xi}_j (x_k,t_k)= \lim_{k\to \infty}
\left[ O(x_k,t_k)\Xi(x_k,t_k)\right]_j=\lim_{k\to \infty}
e_j(x_k,t_k)\cdot \Xi (x_k,t_k)=0
$$
for all $j\in \{1, \cdots, m\}$.
$\square$\\
\ \\
 \noindent{\bf Proof of Theorem 1.2}
 Taking derivative of $\gamma (s,t)=X(\gamma_0 (s ),t)$ with respect
 to $s\in I$, we have
 \bb\label{1.21}
 \frac{\partial \gamma (s,t)}{\partial s} =\nabla_a X(\gamma_0
 (s),t) \frac{\partial \gamma_0 (s)}{\partial s}.
 \ee
Since $W(x,t)$ is transported by $\{ X(\cdot ,t)\}$, we have, along
the curve $t\mapsto\gamma(s,t)$,
 \bb\label{1.22}
 W(\gamma (s,t),t)=\nabla_a X(\gamma _0 (s),t)W_0 (\gamma_0 (s)).
 \ee
 By hypothesis, since $\gamma_0 (s)$ is an integral curve of $W_0
 (\gamma_0 (s))$, there exists $f(s)\neq0$ for all $s\in I$ such that
 \bb\label{1.23}
 \frac{\partial \gamma_0 (s)}{\partial s}=f(s)W_0 (\gamma_0 (s)),
 \ee
 and from (\ref{1.21}) we have
 \bb\label{1.24}
 \frac{\partial \gamma (s,t)}{\partial s}=f(s)\nabla_a X(\gamma _0 (s),t)W_0 (\gamma_0 (s))= f(s)W(\gamma
 (s,t),t),
 \ee
 which shows that $s\mapsto \gamma(s,t)$ is an integral curve  of
 $W(\gamma(s,t) ,t)$ for each $t\in [0, T)$.
From (\ref{1.23}) and (\ref{1.24}) we obtain
 \bb
\frac{1}{|f(s)|}= \frac{|W(\gamma
(s,t),t)|}{\left|\frac{\partial\gamma}{\partial
\sigma}(s,t)\right|}= \frac{|W_0(\gamma_0
(s)|}{\left|\frac{\partial\gamma_0}{\partial \sigma}(s)\right|}.
 \ee
$\square$ \\

\section{Applications to inviscid hydrodynamics}
\setcounter{equation}{0}

We discuss the implications of the previous general theorems on some
of the ideal fluid mechanics equations.

\subsection{The surface quasi-geostrophic equation}

In this subsection  we are concerned on the the following 2D
quasi-geostrophic equation in $\Bbb R^2$.
\[
\mathrm{ (QG)}
 \left\{ \aligned
 &\frac{\partial\theta}{\partial t} +(v\cdot \nabla )\theta =0,
  \\
 & v =-\nabla ^\bot (-\Delta )^{-\frac12} \theta ,\\
  &\theta(x,0)=\theta_0 (x),
  \endaligned
  \right.
  \]
where $\theta=\theta(x_1, x_2, t)$ denotes the scalar temperature,
and $v=(v_1,v_2), v_j=v_j(x_1, x_2,t), j=1,2,$ is the velocity of
the fluid, and $\nabla^\bot =(-\partial_{x_2},
\partial_{x_1})$. Thanks to the pioneering work by Constantin, Majda
and Tabak(\cite{con5}), in particular the observation of its
resemblance to the 3D Euler equations, there are many studies on
(QG)(see e.g. \cite{con1, cor1, cor2, cor4, wu1, wu2} and references
therein).
 Let $\{ X(\cdot,t)\}$ be the particle trajectory mapping generated by
 $v(x,t)$.
 Taking operation of $\nabla^\bot$ on the first equation of (QG) we
 obtain
 \bb\label{2.1}
 \frac{\partial}{\partial t} \nabla^\bot \theta +(v\cdot \nabla
 )\nabla ^\bot \theta =(\nabla^\bot \theta \cdot \nabla )v,
 \ee
from which  we have the transport formula for $\nabla^\bot (x,t)$,
 \bb\label{2.2}
 \nabla^\bot \theta (X(a,t),t)=\nabla_a X(a,t) \nabla^\bot \theta_0
 (a).
  \ee
 As in the
previous section we set the back-to-label map, $A(\cdot,t)=X^{-1}
(\cdot ,t)$ below. The following theorem is immediate from
(\ref{2.2}) and Theorem 1.1, and the fact that det $(\nabla_a
X(a,t))\equiv 1$, which is equivalent to the incompressibility
condition,  div $v=0$.
\begin{thm}
Let $(v(x,t), \theta (x,t))$ be a smooth solution of (QG) with
initial data satisfying $\|\nabla \theta_0 \|_{L^\infty} <\infty$,
which generates the particle trajectory map $\{X(\cdot ,t)\}$ and
the particle trajectory map $A(\cdot, t)$. We set the direction
vector field $ \xi (x,t)= \frac{\nabla^{\bot}\theta
(x,t)}{|\nabla^{\bot}\theta
  (x,t)|},
  $ and let $\{e_1(x,t), e_2(x,t)\}$ and $\{ \lambda_1 (x,t), \lambda_2 (x,t)\}$
   be the normalized eigenvectors and the corresponding eigenvalues
of the matrix
$$ M(x,t)= \left(\nabla A (x,t)\right)^{T} \nabla
A(x,t).$$
 We keep the order of magnitude so that
 $$ \lambda_1 (x,t)>\lambda_2 (x,t) >0 \quad \forall (x,t).$$
 Suppose there exists a sequence $\{(x_k, t_k)\}$ tending
to $(\bar{x}, \bar{t} )$ as $k\to \infty$ such that $ \lim_{k\to
\infty}|\nabla \theta (x_k ,t_k)|= \infty$, then necessarily
 \bb
\lim_{k\to \infty} \lambda_1 (x_k, t_k)=\infty, \quad \lim_{k\to
\infty} \lambda_2 (x_k, t_k)=0,
 \ee
 and
 \bb\label{qg1}
\lim_{k\to \infty} | \xi (x_k, t_k)-e_2 (x_k, t_k) |=0.
 \ee
\end{thm}
We just note that (\ref{qg1}) follows from
 \bb\label{qg2}
  \lim_{k\to \infty}
\xi(x_k, t_k)\cdot e_1 (x_k, t_k) =0
 \ee together with $ (\xi\cdot e_1 )^2 +
(\xi\cdot e_2 )^2 =|\xi |^2=1.$ On the other hand hand, (\ref{qg1})
implies that the direction field tends to align with the direction
of the infinite stretching rate near the possible singularity, while
(\ref{qg2}) shows that the direction of zero compression rate is
orthogonal to it.\\
\ \\
Since any smooth level curve of $\theta_0$ is an integral curve of
$\nabla^\bot \theta_0 $, applying Theorem 1.2 to (QG), we obtain the
following theorem.
\begin{thm}
Let $(\theta (x,t). v(x,t))$ be a smooth solution of (QG), and $\{
X(\cdot, t)\}$ the particle trajectory generated by $v(x,t)$.  Let
$\{ \gamma_0(s)\}_{s\in I}$ be a level curve of $ \theta_0$, which
is also a level curve of $\theta_0$. We set $\gamma (s,t)=X(\gamma_0
(s),t)$, then $\gamma(s,t)$ is also a level curve of $\theta(x,t)$.
Moreover, we have the following invariants along the trajectories of
level curves of $\theta(x,t)$:
 \bb
 \frac{|\nabla^\bot
 \theta (\gamma (s,t),t)|}{ \left|\frac{\partial \gamma (s,t)}{\partial s}
 \right|}= \frac{|\nabla^\bot \theta_0 (\gamma_0(s))|}{
 \left|\frac{\partial \gamma_0(s)}{\partial s}
 \right|}.
\ee
\end{thm}
\begin{cor}
Suppose there exist a sequence $\{ (s_k, t_k )\}$ and  $(\bar{s},
\bar{t})$ such that $(s_k, t_k ) \to (\bar{s}, \bar{t})$, and
 \bb
\lim_{k\to \infty} |\nabla ^\bot \theta(\gamma (s_k, t_k ),t_k )|
=\infty ,
 \ee
  then necessarily
 \bb
\lim_{k\to \infty} \left|\frac{\partial \gamma}{\partial s} (s_k,
t_k )\right| =\infty .
 \ee
Namely, the blow-up of $|\nabla^\bot \theta|$  at a point is
accompanied by an infinite stretching of level curves at the same
point in the tangential direction to the curve.
\end{cor}

\subsection{The Euler equations for isentropic flows}

We are concerned here with the following Euler equations for the
isentropic  fluid flows in $\Bbb R^n$, $n=2,3$,
 \[
 (E) \left\{ \aligned
 &\rho\frac{\partial v}{\partial t} +\rho(v\cdot \nabla )v =-\nabla p ,
  \\
 &\frac{\partial \rho}{\partial t} +\mathrm{div}(\rho v)=0 , \\
 & v(x,0)=v_0 (x), \quad\rho (x,0)=\rho_0 (x),
  \endaligned
  \right.
  \]
where $v=(v_1, \cdots , v_n )$, $v_j =v_j (x, t)$, $j=1, \cdots, n,$
is the velocity of the flow,  $\rho=\rho(x,t)$ is the mass density
of the fluid, $p=p(x,t)$ is the scalar pressure, and $v_0, \rho_0 $
are the given initial velocity and density. The case of homogeneous
incompressible Euler equations corresponds to $\rho(x,t)\equiv
$const., for which we denote by $(E)_0$.  The problems of finite
time blow-up/global regularity for the  systems $(E)$ and $(E)_0$
are both outstanding open problems in the mathematical fluid
mechanics. For $(E)_0$ there are results on the blow-up criterion
initiated by Beale, Kato and Majda(\cite{bea}), and refined by
authors in \cite{cha1,con1, con3, koz1, koz2, maj1}). The study of
the Euler system in terms of the volume preserving maps
 are previously done by  many authors(see e.g. \cite{arn, bre}).
 The geometric type of approaches, emphasizing the role of the direction of
 vorticity  for the regularity/singularity of solutions are studies in
 \cite{con4,den1, den2,gib, cha3},
 the spectral dynamics type of approaches are studied in \cite{liu, cha2}, and
some of the plausible scenarios leading to singularities are
 excluded in \cite{cor3, cor4, cha4, cha5}.
 Let $\{ X(\cdot,t)\}$ be the particle trajectory mapping generated by
 $v(x,t)$, defined by a smooth solution of the solutions of (E), and
 $A(x,t)=X^{-1} (x ,t)$
 be the back-to-label
 map. Let $\o(x,t)=$curl $v(x,t)$ be the vorticity.
 The following vorticity transport formula is well-known(see e.g.
 \cite{cho}) for $(E)$.
 \bb\label{eulertrans}
 \frac{\o (X(a,t),t)}{\rho (X(a,t),t)} =\nabla_a X(a,t) \frac{\o_0 (a)}{\rho_0 (a)}.
 \ee
 Applying Theorem 1.1 to the case of $(E)_0$, for which
  we have
 \bb
 \o (X(a,t),t) =\nabla_a X(a,t) \o_0 (a),
 \ee
 as well as det $(\nabla_a X(a,t))\equiv 1$, we obtain the following
 theorem.
\begin{thm}
Let $\o(x,t)$ be the vorticity of a smooth solution $v(x,t)$ of
$(E)_0$  in $\Bbb R^3$ with initial vorticity satisfying $\|\o_0
\|_{L^\infty} <\infty$.  The particle trajectory map $\{X(\cdot
,t)\}$ and the particle trajectory map $A(\cdot, t)$ are generated
by $v(x,t)$. We set the vorticity direction field $ \xi (x,t)=
\frac{\o(x,t)}{|\o
  (x,t)|}
  $. Let $\{(\lambda_j (x,t), e_j(x,t))\}_{j=1}^3$ be the pairs of
 the eigenvalues and  normalized eigenvectors of the matrix
$$ M(x,t)= \left(\nabla A (x,t)\right)^{T} \nabla
A(x,t),
 $$
  where we keep the ordering for the corresponding eigenvalues
 \bb
  \lambda_1  (x,t)\geq \lambda _2 (x,t)\geq \lambda_3 (x,t)>0 .
  \ee
 Suppose there exists a sequence $\{(x_k, t_k)\}$ tending
to $(\bar{x}, \bar{t} )$ as $k\to \infty$ such that $ \lim_{k\to
\infty}|\o (x_k ,t_k)|= \infty$, then necessarily
 \bb
\lim_{k\to \infty}\lambda_1 (x_k, t_k) = \infty \quad
\mbox{and}\quad \lim_{k\to \infty}\lambda_3 (x_k, t_k) =0,
 \ee
 and
 \bb\label{2.11}
 \lim_{k\to \infty}\xi(x_k, t_k ) \cdot e_1 (x_k, t_k )=0.
 \ee
Furthermore, if
 \bb
 \lim\inf_{k\to \infty} \lambda_2 (x_k,t_k )>0,
\ee
 then
  \bb\label{2.13}
 \lim_{k\to \infty}  \xi (x_k, t_k)\cdot e_2 (x_k,t_k)=0.
 \ee
\end{thm}
\noindent{\textsc{Remark} 2.1} In the case when (\ref{2.11}) and
(\ref{2.13}) happen, we note that
 \bb
  \lim_{k\to \infty} \xi (x_k, t_k)\cdot e_3 (x_k,t_k)=1,
 \ee
 which is equivalent to
 \bb
 \lim_{k\to \infty} | \xi (x_k, t_k)-e_3 (x_k, t_k) |=0.
 \ee
 Namely, as $(x_k, t_k)$ tends to $(\bar{x}, \bar{t})$,
 the sequence of vorticity direction vectors $\{\xi (x_k,t_k)\}$ tends to align
with the eigenvector of $M(x_k,t_k)$ with smallest eigenvalue, which
is in the direction of maximum stretching rate. Taking into account
of the formula (\ref{eulertrans}), we obtain the following theorem
immediately from Theorem 1.2.
\begin{thm}
Let $ (v(x,t), \rho(x,t))$ be a smooth solution of (E), and $\{
X(\cdot, t)\}$ be the particle trajectory generated by $v(x,t)$. Let
$\gamma_0(s)$ be a vortex line for the initial vorticity $\o_0$. We
set $\gamma (s,t)=X(\gamma_0 (s),t)$, which is also a vortex line by
the Helholtz theorem.  Then, we have the following invariants along
the trajectories of the vortex lines:
 \bb
 \frac{|\o  (\gamma (s, t),t)|}{ |\rho (\gamma (s, t),t)|\left|\frac{\partial \gamma (s, t)}{\partial s}
 \right|}= \frac{|\o_0 (\gamma_0 (s))|}{
|\rho_0 (\gamma_0(s))| \left|\frac{\partial \gamma_0(s)}{\partial s}
 \right|}.
\ee
\end{thm}
\begin{cor} Let $\o=\mathrm{curl}\, v$ and $\gamma(s,t)$ as in Theorem 2.4,
and $\left\|\o_0/\rho_0 \right\|_{L^\infty} <\infty$. Suppose there
exist a sequence $\{ (s_k, t_k )\}$ and  $(\bar{s}, \bar{t})$  such
that $(s_k, t_k ) \to (\bar{s},\bar{ t})$, and
 \bb
\lim_{k\to \infty} \frac{|\o (\gamma (s_k, t_k ),t_k
)|}{|\rho(\gamma (s_k, t_k ),t_k)|} =\infty ,
 \ee
then necessarily
 \bb
\lim_{k\to \infty} \left|\frac{\partial \gamma}{\partial s} (s_k,
t_k )\right| =\infty .
 \ee
Namely, a singularity of $|\o|/|\rho|$ at a point is accompanied by
 infinite stretching of the vortex line at the same point.
\end{cor}
\noindent{ \textsc{Remark} 2.2} In the case of $(E)_0$ thin vortex
tube(vortex filament) stretching near the singularity of vorticity
in the 3D Euler equations is well-known fact in the elementary fluid
mechanics(see e.g.\cite{cho}), which is immediate from  Kelvin's
circulation theorem.
The above corollary, on the contrary, is about stretching
of individual vortex lines, not the tubes.

\subsection{The magnetohydrodynamic equations}

We are  concerned here on the the ideal MHD system in $\Bbb R^3$.

\[
\mathrm{ (MHD)}
 \left\{ \aligned
 &\frac{\partial v}{\partial t} +(v\cdot \nabla )v =-\nabla p -b\times \mathrm{curl} \, b, \\
 &\frac{\partial b}{\partial t} +(v\cdot \nabla )b =(b \cdot \nabla )
 v,\\
 &\quad \textrm{div }\, v =\textrm{div }\, b= 0 ,\\
  &v(x,0)=v_0 (x), \quad b(x,0)=b_0 (x)
  \endaligned
  \right.
  \]
where $v=(v_1, v_2, v_3 )$, $v_j =v_j (x, t)$, $j=1,2,3$, is the
velocity of the flow, $p=p(x,t)$ is the scalar pressure, $b=(b_1,
b_2, b_3 )$, $b_j =b_j (x, t)$, is the magnetic field, and
 $v_0$, $b_0$ are the given initial velocity and magnetic field,
 satisfying div $v_0 =\mathrm{div}\, b_0= 0$.
 Below $\{ X(a,t)\}$ is the particle
trajectory mapping generated by
 $v(x,t)$, defined by a smooth solution of (MHD), and $A(x,t)=X^{-1} (x ,t)$ is the back-to-label
 map.
As for the vorticity transport formula the  second equation of (MHD)
implies that we have
 \bb
 b (X(a,t),t) =\nabla_a X(a,t) b_0 (a),
 \ee
 which provides us the following theorem due to Theorem 1.2.
\begin{thm}
Let  $v(x,t)$ be a smooth solution of (MHD) with the initial data
satisfying $\|b_0 \|_{L^\infty} <\infty$, which generates the
particle trajectory map $\{X(\cdot ,t)\}$ and the particle
trajectory map $A(\cdot, t)$.
 We set the direction vector field of
the magnetic field $ \xi (x,t)= \frac{b(x,t)}{|b
  (x,t)|},
  $ and let $\{(\lambda_j (x,t), e_j(x,t))\}_{j=1}^3$ be the pairs of
the eigenvalues  and normalized eigenvectors of the matrix
$$
 M(x,t)= \left(\nabla A (x,t)\right)^{T} \nabla A(x,t),
 $$
  where we keep the ordering for the corresponding eigenvalues
 \bb
  \lambda_1  (x,t)\geq \lambda _2 (x,t)\geq \lambda_3 (x,t)>0 .
  \ee
 Suppose there exists a sequence $\{(x_k, t_k)\}$ tending
to $(\bar{x},\bar{ t})$ as $k\to \infty$ such that $ \lim_{k\to
\infty}|b (x_k ,t_k)|=\infty$, then necessarily
 \bb
\lim_{k\to \infty}\lambda_1 (x_k, t_k) =\infty, \quad
\mbox{and}\quad \lim_{k\to \infty}\lambda_3 (x_k, t_k) =0,
 \ee
 and
 \bb\label{2.22}
 \lim_{k\to \infty}\xi  (x_k, t_k )\cdot e_1 (x_k, t_k )=0.
 \ee
Furthermore, if
 \bb\label{2.23}
 \lim\inf_{k\to \infty} \lambda_2 (x_k,t_k )>0,
\ee
 then
  \bb\label{2.24}
 \lim_{k\to \infty}  \xi (x_k, t_k)\cdot e_2 (x_k,t_k)=0.
 \ee
\end{thm}

\noindent{ \textsc{Remark }2.3} We have similar remark to Remark
2.1, and have
 \bb
 \lim_{k\to \infty} | \xi (x_k, t_k)-e_3 (x_k, t_k) |=0
 \ee
 in case (\ref{2.23}) holds. Namely, as $(x_k, t_k)$ tends to $(\bar{x}, \bar{t})$,
 the sequence of magnetic direction vectors $\{\xi (x_k,t_k)\}$ tends to align
with the eigenvector of $M(x_k,t_k)$ with smallest eigenvalue, which
is in the direction of maximum stretching rate.
\ \\
\begin{thm}
Let $(v(x,t), b(x,t))$ be a smooth solution of (MHD), and $\{
X(\cdot, t)\}$ the particle trajectory generated by $v(x,t)$. Let
$\gamma_0(s)$ be a integral curve  for the initial magnetic field
$b_0 (\gamma_0 (s))$. We set $\gamma (s,t)=X(\gamma_0 (s),t)$, which
is also an integral curve of the magnetic  field $b(x,t)$. Then, we
have the following invariants:
 \bb
 \frac{|b (\gamma (s, t),t)|}{ \left|\frac{\partial \gamma (s, t)}{\partial s}
 \right|}= \frac{|b_0 (\gamma_0 (s))|}{
 \left|\frac{\partial \gamma_0(s)}{\partial s}
 \right|}.
\ee
\end{thm}

\begin{cor}
Let $b, \gamma (s,t)$ be as in Theorem 2.6, and $\|b_0 \|_{L^\infty}
<\infty$. Suppose there exist a sequence $\{ (s_k, t_k )\}$ and
$(\bar{s}, \bar{t})$ such that $(s_k, t_k ) \to (\bar{s},\bar{ t})$,
and
 \bb
\lim_{k\to \infty} |b (\gamma (s_k, t_k ),t_k )| =\infty ,
 \ee
then necessarily
 \bb
\lim_{k\to \infty} \left|\frac{\partial \gamma}{\partial s} (s_k,
t_k )\right| =\infty.
 \ee
Namely, a singularity of the magnetic field of (MHD)  is accompanied
by an infinite stretching of magnetic field lines in the direction
of magnetic field.
\end{cor}

  \end{document}